\documentclass [12 pt]{article}
\usepackage{amsfonts}
\title{Probabilistic Extensions of the Erd\H os-Ko-Rado Property}
\author{Anna Celaya \\
Department of Mathematics\\
University of Wisconsin \\
{\tt accelaya@yahoo.com}\\
\and
Anant P.~Godbole \\
Department of Mathematics\\
East Tennessee State University \\
{\tt godbolea@mail.etsu.edu}\and
Mandy Rae Schleifer\\
Department of Psychology\\
Duquesne University\\
{\tt mandyrae123@yahoo.com}}
\begin{document}
\def\p{\mathbb P}
\def\e{{\mathbb E}}
\def\lr{\left(}
\def\rr{\right)}
\def\var{{\rm Var}}
\def\lc{\left\{}
\def\rc{\right\}}
\def\qed{\vbox{\hrule\hbox{\vrule\kern3pt\vbox{\kern6pt}\kern3pt\vrule}\hrule}}
\def\pnd{{\frac{1}{2}{n\choose k}{n-k\choose k}}}
\def\nk{{n\choose k}}
\def\nkk{{n-k\choose k}}
\def\kr{{k\choose r}}
\def\nkr{{n-k\choose k-r}}
\def\nkj{{n-k\choose k-j}}
\def\hmax{{{k\choose \frac{k^2}{n}} {n-k\choose {k - \frac{k^2}{n}}}}}
\def\kj{{k\choose j}}
\def\pndr{{\frac{1}{2}{n\choose k}{k\choose r}{n-k\choose k-r}}} 
\def\thresh{{\sqrt{\frac{2 \nk}{{k\choose r}{n-k\choose k-r}}}}}
\def\thresha{{\sqrt{\frac{2 \nk}{\nkk}}}}
\def\threshb{{\sqrt{2} e^{\frac{k^2}{2 n}}}}
\newtheorem{thm}{Theorem}
\newtheorem{result}{Result}
\maketitle
\begin{abstract}
The classical Erd\H os-Ko-Rado (EKR) Theorem states that if we choose a family of subsets, each of size \(k\), from a fixed set of size \(n\ (n > 2k)\), then the largest possible pairwise intersecting family has size \(t ={n-1\choose k-1}\).  We consider the probability that a randomly selected family of size \(t=t_n\) has the EKR property (pairwise nonempty intersection) as $n$ and $k=k_n$ tend to infinity, the latter at a specific rate.  As $t$ gets large, the EKR property is less likely to occur, while as $t$ gets smaller, the EKR property is satisfied with high probability.  We derive the threshold value for $t$ using Janson's inequality.
Using the Stein-Chen method we show that the distribution of $X_0$, defined as the number of disjoint pairs of subsets in our family, can be approximated by a Poisson distribution.  We extend our results to yield similar conclusions for $X_i$, the number of pairs of subsets that overlap in exactly $i$ elements.  Finally, we show that the joint distribution $(X_0, X_1, \ldots, X_b)$ can be approximated by a multidimensional Poisson vector with independent components.
\end{abstract}

\section{Introduction}
The classical combinatorics literature is replete with fundamental results on properties of intersecting families of sets and subsets of fixed element sets.  Results in this genre include Sperner's theorem \cite{wilson}, Kneser's theorem \cite{greene} and the starting point of this paper, the Erd\H os-Ko-Rado theorem.  In 1960, Erd\H os and Rado proved that for each pair of positive integers $n$ and $k$, with $k \geq 2$, there corresponds a least positive integer ${\cal {F}}(n,k)$ such that if ${\cal {F}}$ is a family of more than ${\cal {F}}(n,k)$ sets, each set with $n$ elements, then some $k
$ of the sets have pairwise the same intersection.  Together with Ko in 1961, they produced the Erd\H os-Ko-Rado (EKR) Theorem \cite{ekr}, which states that if ${\cal {F}}$ is a pairwise intersecting family of $k$-element subsets chosen from an $n$-element set with $n \geq 2k$ then $\|{\cal {F}}\| \leq {n-1\choose k-1}$, where, throughout this paper we denote the cardinality of a set $A$ by $\|A\|$.  Note that  a maximal ensemble of this type may be constructed by selecting all $k$-subsets that contain a fixed element $a$.

In this paper, a family of sets defined to have the ``EKR property" is a group of $k$-sized subsets chosen from an $n$-element set such that there exists a nonempty intersection between any two subsets.  We address the following threshold-type question: as $n$ and $k$ tend to infinity, how many $k$-subsets will we be able to choose (at random, using the uniform probability measure on $k$-sets) such that the EKR property is ``almost always" or ``almost never" satisfied, where these terms are used in the graph-theoretic rather than measure theoretic sense?     

Let $R$ denote our family of $k$-sets, so that $\|R\| = t$ is the size of our family of $k$-sets.  We will let $X_0$ be the random variable representing the number of disjoint pairs in our selection of $t$ subsets.  $X_0 = 0$ corresponds to no disjoint pairs being in our family, i.e., to this family of sets having the EKR property.  We will, accordingly, often use the notation $\p({\rm EKR}) $ instead of $  \p(X_0 = 0).$  

In Section 2, we use Janson's exponential inequalities, as found, e.g., in \cite{alon}, to prove asymptotic threshold results of the form 

\noindent $\p({\rm EKR})\to1\enspace (t\ll t_0)$ and $\p({\rm EKR}) \to0\enspace(t\gg t_0)$, where, throughout this paper, we write $a_n\ll b_n$ or $b_n\gg a_n$ if $a_n/b_n\to0\enspace(n\to\infty)$.  In Section 3, we will examine the distribution of $X_0$ and use the Stein-Chen method of Poisson approximation \cite{bhj} to prove that $d_{TV}({\cal L}(X_0), {\rm Po}(\e(X_0))\to0$ as $n,k=k_n\to\infty$ at an appropriate rate, where ${\cal L}(Z)$ represents the distribution of the random variable $Z$, $d_{TV}$ the usual total variation distance, and ${\rm Po}(\lambda)$ the Poisson distribution with parameter $\lambda$.  The generalization of the EKR property alluded to in the abstract will be provided in Section 4, where we present asymptotic results on the existence and numbers of pairs of $k$-sets that overlap in exactly $r$ elements.  Finally, in Section 5, we discuss the joint distribution of the ensemble ${(X_0, X_2, . . . , X_b)}$ where for $i = 0, 2, \ldots, b$, $X_i$ is the random variable representing the number of pairs of $k$-sets which overlap in exactly $i$ elements.

We end this section with a few potential applications:  Suppose that at an international event, there is a need for interpreters to be hired.  If $n$ is the total number of languages spoken amongst the $t$ interpreters who will be present at such an event and if each interpreter speaks $k$ languages, the results of this paper could be used to determine thresholds for $t$ so that any two interpreters can converse with each other.  Similarly, we may consider a workshop with $t$ participants, each of whom is randomly scheduled to attend $k$ sessions out of a total of $n$.  The EKR property would suggest that any pair of participants could have a meaningful dinner conversation, while the results of this paper would enable one to derive probabilistic conclusions along the same lines, such as the following:  What value of $t$ would ensure with probability at least 0.95 that at least $x$ pairs of participants can only talk about a single session that they both attended, and that between $y$ and $z$ pairs of participants find that they attended between two and five common sessions?

\section{Threshold for the EKR Property}
Intuitively, one would imagine that for appropriately chosen values of $n$ and $k$, a small number of randomly chosen $k$-sets would allow the EKR property to hold with high probability, whereas even a ``slightly" larger collection  would cause the pairwise intersection property to be ruined.  We make this precise in the following result, which makes use of Janson's exponential inequality \cite{alon}.

\begin{thm}

Let $t$ denote the number of $k$-sets chosen at random from an $n$-element set. We set 
$$X_0=\sum_{j=1}^\pnd I_j,$$ where $I_j$ equals 1 if the $j$th pairwise disjoint pair is in the selected ensemble and $I_j=0$ otherwise. Then with $t_0={\sqrt{{2\nk}/{\nkk}}}$, 
$$\ \ \ \mbox{P($X_0 = 0$)} \to\left\{ \begin{array}{ll}
	1 & \mbox{if $t \ll t_0$} \\
	0 & \mbox{if $t \gg t_0$} \\
	e^{-A^2} & \mbox{if $t = (A+o(1))t_0$}
	\end{array}
\right. $$
as $n,k \to \infty$, provided $k \gg \sqrt{n}, k=o(n).$  If $k=o(n^{2/3})$, we may use the more convenient $t_0 = \sqrt{2}e^{\frac{k^{2}}{2n}}$ in the above result.\par 

\end{thm}

\noindent{\bf Proof}  We start by altering our model slightly and choosing each $k$ set independently with probability $p=t/\nk$. In other words, we flip a coin with bias $p= {t}/{{n\choose k}}$, to decide whether each of the ${n\choose k}$ subsets will be in our family of $k$-sets.  We thus obtain a random collection $R$ of $k$-sized subsets where $\e(\| R \|) = t$. Janson's inequality, which bounds the probability that none of a sequence of ``undesirable" events $B_i; i\in I$ occurs, asserts that under certain fairly general conditions

$$ \prod_{i\in I} \p (\overline{B}_i) \leq \p\left( \bigcap_{i\in I} \overline{B}_i \right) \leq \exp \left( - \mu + \frac{1}{1-\epsilon} \frac{ \Delta}{2} \right)
$$
where $ {\p(B_i)} \leq \epsilon$ for each $i$,  $ \mu = \sum_{i \in I}P(B_i) $,  and, for $\sim$ to be defined below, 
$$ \Delta = \sum_{i \sim j} \p(B_i \cap B_j). $$

Let $B_i$ be the event that the $i$th disjoint pair of subsets is in our selection of $k$-sets.  It follows that $\p\left( \bigcap_{i\in I} \overline{B}_i \right)$ = $\p$(Selecting a family with the EKR property) = $\p(X_0 = 0)$.  
Following the canonical set-up for the validity of the Janson inequality, we say $i \sim j$ if the $i$th and $j$th disjoint pairs of $k$-sets have one set in common.  The probability $\p(B_i$) that a particular pair of disjoint sets is in our selection of $k$-sets is $p^2$.  Since there are $\pnd$ possible disjoint pairs, it follows that $\mu=\e(X_0)=\pnd p^2$ and that
$$\prod_{i = 1}^{\pnd} {\p} (\overline{B}_i) 
= (1-p^2)^\pnd.$$
Using the inequality $1-x\ge\exp{\lc-\frac{x}{1-x}\rc}$, the lower Janson inequality yields
$$  \exp\lc- \frac{\pnd p^2}{1-p^2}\rc \leq \p \left( \bigcap_{i\in I} \overline{B}_i \right),$$
so that $\p({\rm EKR})\to1$ when $p = o\left(\sqrt{{2}/{\nk\nkk}}\right)$ or equivalently when $t=o\lr{\sqrt{{2\nk}/{\nkk}}}\rr$.  Note:  There is a simpler proof of this fact using Markov's inequality, but we have presented the above proof for uniformity of exposition. Also, we have assumed above, as we will throughout this paper, that $p\to0$. 

In terms of an equivalent (and more convenient) exponential bound, we have, since $\sqrt{2e^{\frac{k^2}{n}}} \leq \sqrt{{2 \nk}/{\nkk}}$, 
that
$$\p({\rm EKR})\to1\enspace{\rm when}\enspace t\ll \sqrt{2e^{\frac{k^2}{n}}}.$$
Let us now see when when the upper bound in Janson's inequality tends to 0, thus yielding $\p({\rm EKR})\to0$.  We have 
$$\Delta= \sum_{i \sim j} \p(B_i \cap B_j)\le\pnd\cdot2\cdot\lr\nkk-1\rr \cdot p^3\le\nk\nkk^2p^3,$$
so that the upper Janson inequality yields
$$
\p\left( \bigcap_{i\in I} \overline{B}_i \right) \leq  \exp \left(\frac{-\pnd p^2 (1 - p^2 - \nkk p)}{1 - p^2} \right).
$$
We note that this upper bound tends to zero when $t = p \nk$ satisfies
$$ \thresha \ll t \le \frac{ \nk}{2\nkk}.$$
Assuming, as stated in the theorem, that $k^2 \gg n$ and $k \ll n$, it is easy to verify that the above range for $t$ is a valid one, i.e., that
${\sqrt{{2\nk}/{\nkk}}} \ \mbox{is indeed} \ \ll  { \nk}/{2\nkk}$
We conclude by monotonicity that 
$$ \p(X_0 = 0) \to 0 \ \mbox {when}\  t \gg {\sqrt{\frac{2 \nk}{\nkk}}}.$$
It is easy to verify using the inequality $1-x\ge\exp\{-x/(1-x)\}$ that 
$$\frac{\nkk}{\nk}\ge\exp\lc-\frac{k^2}{n}-\frac{2k^3}{n^2}\rc,$$ 
so that $\p({\rm EKR})\to0$ when $ t \gg \threshb $, provided that 
$k=o(n^{2/3})$, as asserted.

Next, we examine the behavior of the EKR property around this threshold value. If we let $t = (A+o(1)){\sqrt{\nk/\nkk}} $, and let $n$ and $k$ go to $\infty$, the lower and upper bounds of the Janson inequality together yield
$\p(X_0=0)\sim e^{- A^2} $, proving the last part of the theorem, {\it but for the altered model}, i.e. when we {\it expect} to chose $t=p\cdot\nk$ subsets.  It remains to ``derandomize" our results, so as to verify that the same threshold is valid when {\it exactly} $t$ subsets are chosen.  We proceed in a fashion similar to that in \cite{janson}.

First we derandomize the result corresponding to $P(X_0 = 0) \to 0$.  Let $\|R\|$ be the exact size of the chosen family.  Assuming that $p \gg \sqrt{{2}/{\nk \nkk}}$, we wish to prove that if $\|R\| = p \nk$ $k$-sets are chosen, then $\p(X_0 = 0) \to 0$.  We have by monotonicity

\begin{eqnarray*}
\p\left( X_0 = 0  \bigg\vert \|R\| = p \nk \right) &\leq& \p\left( X_0 = 0 \bigg\vert  \|R\| \leq p \nk \right)\\
&\leq& \frac{\p (X_0 = 0)}{\p(\|R\| \leq p \nk)}\\
&\le&3\p(X_0=0)\to0,
\end{eqnarray*}
by our preliminary result, the fact that $\p(A\vert B)\le \p(A)/\p(B)$ and the fact that the central limit theorem (or the approximate and asymptotic equality of the mean and median of the binomial distribution ${\rm Bin}(n,p)$) implies that $\p({\rm Bin}(n,p)\le np)\ge1/3.$  It follows that
$$ \p \left( X_0 = 0 \bigg\vert \|R\| = p \nk \right) \to 0 \ \mbox{if} \ p  \gg \sqrt{\frac{2}{\nk \nkk}},$$
\noindent or
$$ \p\left( X_0 = 0 \bigg\vert \|R\| = t, \ t \gg \thresha \right) \to 0.$$
\noindent To derandomize the case where $\p(X_0 = 0) \to 1$, we proceed as before:
\begin{eqnarray*}
\p\left( X_0 \ge1 \bigg\vert  \|R\| \ = p \nk \right) &\leq& \p \left( X_0 \ge1 \bigg\vert  \|R\| \ \geq p \nk \right)\\
&\leq& \frac{\p (X_0 \ge1)}{\p (\|R\| \geq p \nk)}\\
&\le&3\p(X_0\ge1)\to0
\end{eqnarray*}
if $p\ll{\sqrt{{2}/{\nk\nkk}}}$.  It follows that
$$ \p\left( X_0 = 0 \bigg\vert \|R\| = t, \ t \ll \thresha \right) \to 1,$$
as required.  

Finally, we  know that when  $p = (A+o(1)) \sqrt{{2}/{\nk \nkk}}$ (where $A>0$ is a constant),  $\p(X_0 =0) \to e^{- A^2}$.  We define, with hindsight (but somewhat arbitrarily),

$$ p^{+} = \frac{\left( t + \lc\nk/\nkk\rc^{5/16} \right)}  {\nk}$$ 
and
$$p^{-} = \frac{\left( t - \lc\nk/\nkk\rc^{5/16} \right)}  {\nk}, $$ 
where $t = (A+o(1)) {\sqrt{2\nk/\nkk}}$. Note that both $p^+$ and $p^-$ are of the form  $(A+o(1)) \sqrt{{2}/{\nk \nkk}}$.   
Using first $p^{+}$ as the probability of picking any $k$-set thereby obtaining a random collection $R^+$, we see that 
\[
\e(\| R^+ \|) = p^+\nk=(A+o(1)){\sqrt{2\nk/\nkk}}\]
Note also that

$$ \var(\| R^+ \|) = \e(\| R^+ \|) (1 - p^{+}) \approx \e(\| R^+ \|) \sim (A+o(1)) {\sqrt{2\nk/\nkk}}. $$
It follows that
\begin{eqnarray*}
\p( \| R^+ \| \le t) &=& \p\left( \| R^+ \| \le\e(\|R^+\|) - \lc\nk/\nkk\rc^{5/16} \right) \\
&\leq& \p\left( \bigg\vert \| R^+ \| - \e (\| R^+ \|) \bigg\vert \ge \lc\nk/\nkk\rc^{5/16} \right) \\
& \leq& \frac{(A+o(1)){\sqrt 2}}{\lc\nk/\nkk\rc^{1/8}}\to0,
\end{eqnarray*}
where the final step follows from Chebychev's inequality. 
We thus have
$$\p(X_0 \ge1\bigg\vert\| R^+ \| = t) 
\leq \frac{\p(X_0 \ge1)}{\p(\| R^+\| \geq t)}\to 1-e^{-A^2},$$
implying that 

$$ \liminf_{n\to\infty} \p({\rm EKR} \bigg\vert\| R \| = t) \ge e^{- A^2} $$
The proof using $p^{-}$ follows a similar path, yielding 
$$\limsup_{n\to\infty} \p({\rm EKR} \bigg\vert\| R \| = t) \le e^{- A^2}$$
and thus that
$$ \p\left( {\rm EKR} \bigg\vert\|R\| = t, \ t = (A+o(1)) {\sqrt{2\nk/\nkk}}\right) \to e^{- A^2}.$$

\noindent With all three preliminary results derandomized, the proof of Theorem 1 is complete.\hfill\qed

To provide a numerical comparison, we note that when $k=n^{3/5}$, the Erd\H os-Ko-Rado theorem yields $$t\ge{{n-1}\choose{k-1}}\sim\frac{\exp\{(2n^{3/5}\log n)/5(1+o(1))\}}{{\sqrt{2\pi}}n^{7/10}}\Rightarrow\p(X_0=0)=0,$$
whereas Theorem 1 yields a threshold at ${\sqrt2}e^{{k^2}/2n}\sim{\sqrt2}e^{{n^{1/5}}/2}$.

\section{The Distribution of $X_0$}

In this section, we use the Stein-Chen method of Poisson Approximation \cite{bhj} to prove the following theorem.

\begin{thm}

Consider a family of $k$-element subsets of an $n$-element set, obtained by randomly and independently selecting each $k$-set with probability $p$.  Let $X_0$ represent as before the number of disjoint pairs in our selection of subsets.  Then the  distribution of $X_0$ can be closely approximated by a Poisson distribution with parameter $\lambda = \pnd p^2$ when $p = o {\left( \frac{1}{\nkk} \right)}$.

\end{thm}

\noindent{\bf Proof.}  One of several  approximation theorems in \cite{bhj} (Corollary 2.C.4) states that if we consider a sum $Z = \sum_{j\in I} I_j$ of indicator random variables with $E(Z) = \lambda$, and if for each $j$ there exists a sequence of indicator variables, $J_{ij}$, such that 

\begin{equation}{\cal {L}}(J_{ij} : i \in I) = {\cal{L}}(I_{i} : i \in I \bigg\vert I_j = 1),\end{equation}
and such that for all $i \neq j$, $J_{ij} \geq I_i$, then
$$ d_{TV} \left({\cal {L}}\left( Z \right),{\rm Po} \left( \lambda \right) \right) \ \leq \ \frac{1-e^{-\lambda}}{\lambda} \left( \var\left(Z \right) - \lambda + 2 \sum {\p}^2 \left(I_{j} = 1 \right) \right),$$  
where $d_{TV}$ represents the usual total variation distance and ${\rm Po}(\lambda)$ the Poisson random variable with mean $\lambda$.  In other words, if a coupling exists such that the indicator random variables $I_i$ and $J_{ij}$ are positively related, then the total variation distance between the distribution of the random variable $Z$ and a Poisson distribution with parameter $\lambda$ may be bounded solely in terms of the first two moments of $Z$.

For our problem, we employ the following coupling that clearly satisfies (1):
If $I_j = 1$, we let $J_{ij} = I_i \ \  \forall i$.
If $I_j = 0$, we add one or both unchosen $k$-sets to our collection by changing the associated coin flips as needed.  Then, we set $J_{ij} = 1$ if the addition of these $k$-sets creates the selection of the $i$th pair of disjoint $k$-sets to our collection of sets; it is obvious that 
 $J_{ij} \geq I_i \ \ \forall i \neq j$. Therefore, the above bound may be applied with $\lambda = \e(X)=\pnd p^2$ to yield

\begin{eqnarray*}
d_{TV} \left({\cal {L}}(X_0),{\rm Po}(\lambda) \right) 
&\leq& \frac{\var(X_0)}{\lambda} \ \  -1 + \frac{2 \sum \p^2(I_j=1)}{\lambda}
\\&=& \frac{\var(X_0)}{\lambda}- 1 + 2p^2.
\end{eqnarray*}
Since
\begin{eqnarray*}
\var\left( X_0 \right) &=& \var\left( \sum I_j \right)
\\&=& \sum \var\left( I_j \right) \ \ + \ \   2 \sum _{i<j}{\rm Cov}\left( I_i I_j \right)
\\&\le& \pnd \left( p^2 - p^4 \right) \ \ + \ \ 2 \left(\frac{1}{2} \right) \nk {\nkk}^2 \left( p^3 - p^4 \right),
\end{eqnarray*}
it follows that
$$d_{TV} \left({\cal {L}}(X_0),{\rm Po}(\lambda) \right) \leq 2 \nkk p\ \ +\ \ \left( 1\ -\ 2\nkk \right)p^2,
$$
establishing the result.
\hfill\qed

\noindent Note that the threshold value for $p$ in Theorem 1 does in fact fall within the domain of applicability of Theorem 2.
\section{Pairwise $r$-Overlapping Sets}
We are not motivated, in this section or the next, by combinatorial results such as the Erd\H os-Ko-Rado theorem.  Instead our focus turns to the probabilistic nuances of pairwise intersection properties of randomly selected $k$-sets.  In this section, we use methods similar to those employed in Sections 2 and 3 to prove 

(i) threshold results for the existence of, and 

(ii) distributional results for the numbers $X_r$ of, 

\noindent pairs of $k$-sets that overlap in $r$ elements, $r\ge1$.

\begin{thm}
Let $t$ denote the number of $k$-sized sets chosen at random from an $n$-element set. Let $X_r$ be the number of pairs of sets in the chosen family which overlap in exactly $r$ elements.  Then with $t_0 = {\sqrt{2\nk/\kr\nkr}}$, 
$$\mbox{P($X_r = 0$)} \to\left\{ \begin{array}{ll}
	1 & \mbox{if $t \ll t_0$} \\
	0 & \mbox{if $t \gg t_0$} \\
	e^{-A^2} & \mbox{if $t = (A+o(1))t_0$}
	\end{array}
\right. $$
as $n,k \to \infty$, provided $k \gg \sqrt{n}; k=o(n).$
\end{thm}
\noindent{\bf Proof.}
Again, we use the Janson inequality.  This time, we let $B_i$ be the event that both members of the $i$th pair of subsets which overlap in exactly $r$ elements are among the selected  $t$ subsets. Let $X_r$ be the  number of events $B_i$ which occur, so that $\p\lr\bigcap_{i\in I}\overline{B}_i\rr= \p(X_r = 0)$.
Choose a random collection $R$ of $k$-sized subsets with $\e(\| R \|) = t$ by independently selecting each $k$-set to be in our ensemble with probability $p=t/\nk$. 
We thus have
$$
\prod_{i\in I} \p(\overline{B}_i) 
= (1-p^2)^{\frac{1}{2}{n\choose k}{k\choose r}{n-k\choose k-r}},
$$
so that the lower Janson inequality yields
$$\exp\lc- \frac{\pndr p^2}{1-p^2}\rc \leq  (1-p^2)^{\pndr} \leq \p\left( \bigcap_{i\in I} \overline{B}_i \right),$$
and thus to the conclusion that $\p(X_r=0)\to1$ when $p = o\left(\sqrt{{2}/{\nk\kr\nkr}}\right)$ or equivalently, when $t = o\left({\sqrt{2\nk/\kr\nkr}}\right)$.

Note next that
$ \mu = E( X_r ) = \sum_{i \in I}P(B_i) = \pndr p^2$ and that $\Delta\le{\nk} {\kr}^2 {\nkr}^2 p^3$, so that
the upper bound of Janson's inequality yields
$$
\p\left( \bigcap_{i\in I} \overline{B}_i \right) \leq 
\exp \left(\frac{-\pndr p^2 (1 - p^2 - \kr \nkr p)}{1 - p^2} \right),
$$
and thus to the conclusion that 
$$ \thresh \ll t \le \frac{ \nk}{2{k\choose r}{n-k\choose k-r}}\Rightarrow \p(X_r=0)\to0.$$
Assuming, as before, that $k^2 \gg n$ and $k \ll n$, in order to prove that the above range for $t$ is a valid one we must prove that 
$$ \thresh \ \mbox{is indeed} \ \ll \frac{ \nk}{2{k\choose r}{n-k\choose k-r}}, $$
or equivalently that
$$ \frac{{k\choose r}{n-k\choose k-r}}{\nk} \to0. $$
Now the above is simply the ``hypergeometric" probability  function (making $k$ without replacement selections from a drawer with $k$ white and $n-k$ black socks).  The probability of drawing $r$ white socks is maximized around the mean value of $r = {k^2}/{n}$, and we thus need to show that $ {{k\choose r}{n-k\choose k-r}} \ll \nk$ when $r=k^2/n$. This is confirmed 
below using Stirling's formula and an auxiliary result on Poisson approximation:  By Theorem 6.A in \cite{bhj}, the total variation distance between the distribution of our hypergeometric random variable $W$ and a Poisson distribution with the same mean $k^2/n$ satisfies
\[d_{TV} \left({\cal {L}}(W),{\rm Po}(k^2/n) \right)\le \frac{3k}{n}, \]
so that
\[{\rm Po}(k^2/n,k^2/n)-3k/n\le\p(W=k^2/n)\le{\rm Po}(k^2/n,k^2/n)+3k/n.\]
Stirling's formula yields ${\rm Po}(k^2/n,k^2/n)\sim \frac{\sqrt{n}}{{\sqrt{2\pi}}k},$ so that for $r=k^2/n$,
$$ \frac{{k\choose r}{n-k\choose k-r}}{\nk} \in \frac{\sqrt{n}}{{\sqrt{2\pi}}k}\pm \frac{3k}{n}$$
implying that 
$$\frac{{k\choose r}{n-k\choose k-r}}{\nk}\to0\enspace (n,k \to \infty).$$
We conclude by monotonicity that 
$$ \p(X_r = 0) \to 0 \ \mbox {when}\  t \gg {\sqrt{\frac{2 \nk}{{k\choose r}{n-k\choose k-r}} }}.$$
Finally, if we let $t = (A+o(1)){\sqrt{2\nk/\kr\nkr}} $, the lower and upper bounds of Janson's inequality together yield
$ \p(X_r=0)\to e^{- A^2}\enspace(n,k\to\infty). $
Derandomization of these preliminary results follows as in the proof of Theorem 1.  We only provide details for the last case, {\it viz.}, when the {\it actual} number of $k$-sets chosen is $t = (A+o(1)){\sqrt{2\nk/\kr\nkr}}$.  We know that when  
\begin{equation}p = (A+o(1)) \sqrt{\frac{2}{\nk \kr \nkr}},\end{equation} $\p(X_r =0) \to e^{- A^2}$.  We define, again quite arbitrarily,
$$ p^{+} = \frac{\left( t + {\left( \frac{\nk}{\kr \nkr} \right)}^{\frac{3}{8}}  \right)}{\nk} $$ 

$$ p^{-} = \frac{\left( t - {\left( \frac{\nk}{\kr \nkr} \right)}^{\frac{3}{8}}  \right)}{\nk} $$ 
where $t =(A+o(1)) \thresh$, and note that both $p^+$ and $p^-$ satisfy (2).
Using first $p^{+}$ to yield a random collection $R^+$, we see that
\begin{eqnarray*}
\e(R^+) = p^{+} \nk &=& t + {\left( \frac{\nk}{\kr \nkr} \right)}^{\frac{3}{8}} \\
&=& (A+o(1)) \thresh  + {\left( \frac{\nk}{\kr \nkr} \right)}^{\frac{3}{8}}  \\
&\sim& {\sqrt{2}}A {\left( \frac{\nk}{\kr \nkr} \right)}^{\frac{1}{2}}.
\end{eqnarray*}
Also,
$$ \var(R^+) = \e(R^+) (1 - p^{+}) \sim {\sqrt{2}}A{\left( \frac{\nk}{\kr \nkr} \right)}^{\frac{1}{2}}, $$
and thus
\begin{eqnarray*}
\p( \| R^+ \| \le t) &=& \p \left( \|R^+\| \le \e \|R^+\| - {\left( \frac{\nk}{\kr \nkr} \right)}^{\frac{3}{8}}  \right) \\
&\leq& \p \left( \mid \|R^+\| - \e \|R^+\| \mid \ge {\left( \frac{\nk}{\kr \nkr} \right)}^{\frac{3}{8}} \right) \\
&\le&\frac{{\sqrt 2}A}{\lr{\frac{\nk}{\kr \nkr}}\rr^{\frac{1}{4}}}\to0,
\end{eqnarray*}
\noindent where the final step is true by Chebychev's inequality.  
We thus have
\[
\p(X_r \ge1 \ \mid \ \| R^+ \| = t) \le
\frac{\p (X_r \ge1)}{\p(\| R^+\| \geq t)}\to1-e^{-A^2},
\]
so that 
$$ \liminf \p(X_r = 0 \ \mid \ \| R \| = t) \ge e^{- A^2}. $$
The proof for $p^{-}$ follows similarly, yielding 
$$\limsup \p (X_r = 0 \ \mid \ \| R \| = t) \le e^{- A^2},$$
and consequently that
$$ \p \left( X_r = 0 \ \mid \ \|R\| = t, \ t = (A+o(1)) \thresh \right) \to e^{- A^2}$$
as required. \hfill\qed

\begin{thm}

Consider a family of $t$ subsets, each of size $k$, taken from an $n$-element set.  Let $X_r$ be the random variable which represents the number of pairs of chosen subsets which overlap in exactly $r$ elements.  Then the distribution of $X_r$ can be closely approximated by a Poisson distribution with $\lambda_r = \e(X_r) = \pndr p^2$ when $\  p = o \left( \frac{1}{{k\choose r}{n-k\choose k-r}} \right)$.

\end{thm}

\noindent{\bf Proof.}
Again we use the Stein-Chen method.
Consider the following coupling:
If $I_j = 1$, i.e., if the $j$th pair of $k$-subsets that overlap in $r$ elements is selected, we let $J_{ij} = I_j \ \ \forall i$.  
If $I_j = 0$, we add one or both unchosen $k$-sets to our collection by changing the coin flips as needed.  Then, we set $J_{ij} = 1$ if the addition of these $k$-sets creates the selection of the $i$th pair of $k$-sets which overlap in $r$ elements to our ensemble. 
Since $J_{ij} \geq I_i \ \ \forall i \neq j$, we may apply the same Stein-Chen approximation theorem as before.  Note that  $\lambda_r = \e(X_r)=\frac{1}{2} \nk \kr \nkr p^2$, so that
\begin{eqnarray}
d_{TV} \left({\cal {L}}(X_r),{\rm Po}(\lambda_r) \right) 
&\leq& \frac{1-e^{-\lambda_r}}{\lambda_r} \left( {\rm Var} \left(X_r \right) - \lambda_r + 2 \sum \p^2 \left(I_{j} = 1 \right) \right)\nonumber\\&\le& \frac{{\rm Var}(X_r)}{\lambda_r}  -1+2p^2.
\end{eqnarray}
As before, we see that
\[
{\rm Var} \left( X_r \right) \le \frac{1}{2} \nk \kr \nkr \left( p^2 - p^4 \right) +2 \left( \frac{1}{2} \right)  \nk \kr^2 \nkr^2 \left( p^3 - p^4 \right)
\]
so that (3) yields
\[
d_{TV} \left({\cal {L}}(X_r),{\rm Po}(\lambda_r) \right) \leq 2 \kr \nkr p\ \ +\ \ \left( 1 - 2\kr \nkr \right)p^2,
\]
as needed.
\hfill\qed

\section {The Joint Distribution of ${X_1, X_2, . . . , X_b}$}
\begin{thm}
The joint distribution of $(X_0, X_2, ... X_b)$, can be approximated by $b$ independent Poisson distributions, provided $b$ is ``not too large" and provided the probability $p$ of choosing any particular $k$-set satisfies 
$$ p = o {\left( \frac{1}{\nk {\left( \sum_{j=0}^b \kj \nkj \right) }^2 } \right) }^{\frac{1}{3}} $$
In other words, the total variation distance $d_{TV}$ between these two distributions satisfies 
$$ d_{TV} ( {\cal L} (X_0, X_1, . . . , X_b), \prod^b_{j=0} {\rm Po}( {\lambda}_j ) ) \leq {\epsilon}_{n,b,k} $$
where $ {\epsilon}_{n,b,k} \to 0 $ as $n,k \to \infty$. 

\end{thm}

\noindent{\bf Proof.}
First we note that the restriction on $p$ satisfies the $p$-requirements for each individual Poisson approximation.  This can be seen by noting that for each $r=0,1,\ldots,b$, 

$$ {\left( \frac{1}{\nk {\left( \sum_{j=0}^b \kj \nkj \right) }^2 } \right) }^{\frac{1}{3}} \le\left( \frac{1}{\kr \nkr} \right).$$ 
Consider the indicator variables $I_{\scriptscriptstyle (j,i)}$ where the $(j,i)^{\rm th}$ pair is the $i^{\rm th}$ pair of subsets which overlap in exactly $j$ elements.  If $I_{\scriptscriptstyle (j,i)} = 1$, which is to say that the $(j,i)^{\rm th}$ pair was chosen, then we let $J_{\scriptscriptstyle (\beta , \alpha )}^{\scriptscriptstyle (j,i)} = I_{\scriptscriptstyle (\beta , \alpha )}$.  If $I_{\scriptscriptstyle (j,i)} = 0$ then we choose the $(j,i)^{\rm th}$ pair of subsets.  If $I_{\scriptscriptstyle (\beta , \alpha )} = 1$ after these additional choices, we set $J_{\scriptscriptstyle (\beta , \alpha )}^{\scriptscriptstyle (j,i)} = 1$.

Since we have found a coupling which satisfies 
$$ {\cal L} ( J_{\scriptscriptstyle (\beta , \alpha )}^{\scriptscriptstyle (j,i)} ) = {\cal L} ( I_{\scriptscriptstyle (\beta , \alpha )} \ \mid \ I_{\scriptscriptstyle (j,i)} = 1), $$
Theorem 10.J of \cite{bhj} yields, with $N_j=\frac{1}{2}{n\choose k}{k\choose j}{{n-k}\choose{k-j}},$ 

$$ d_{TV} ( {\cal L} (X_0, X_1, . . . , X_b), \prod^b_{j=0} {\rm Po}( {\lambda}_j ) ) \leq {\epsilon}_{n,b,k}, $$

\noindent where 

$$ {\epsilon}_{n,b,k} = \sum_{j=0}^b \sum_{i=1}^{N_j} \left\{ \p^2(I_{\scriptscriptstyle (j,i)} = 1) + \p(I_{\scriptscriptstyle (j,i)} = 1) \sum_{{ 0\le\beta \le b;1\le\alpha\le N_\beta}\atop{(\beta,\alpha)\ne(j,i)}} \p( I_{\scriptscriptstyle (\beta , \alpha)} \neq J_{\scriptscriptstyle (\beta , \alpha )}^{\scriptscriptstyle (j,i)} ) \right\}. $$
We thus have
\begin{eqnarray*}
{\epsilon}_{n,b,k} &\le& \sum_{j=0}^b \sum_{i=1}^{N_j} \left\{ p^4 + p^2 \sum_{{ 0\le\beta \le b;1\le\alpha\le N_\beta}\atop{(\beta,\alpha)\ne(j,i)}} \p( I_{\scriptscriptstyle (j,i)} = 0 \ \cap \ I_{\scriptscriptstyle (\beta ,\alpha)} = 0 \ \cap \ J_{\scriptscriptstyle (\beta , \alpha )}^{\scriptscriptstyle (j,i)} = 1) \right\} \\
&\le& \sum_{j=0}^b \sum_{i=1}^{N_j} \left\{ p^4 + p^2 \sum_{{ 0\le\beta \le b;1\le\alpha\le N_\beta}\atop{(\beta,\alpha)\ne(j,i);\vert(\alpha,\beta)\cap(j,i)\vert=1}}  2p \right\} \\
&\le& \sum_{j=0}^b \sum_{i=1}^{N_j} \left\{ p^4 + 4p^3 \sum_{\beta=0}^b {k \choose \beta} {n-k \choose k- \beta} \right\} \\
&=& \left( p^4 + 4p^3 \sum_{\beta=0}^b {k \choose \beta} {n-k \choose k- \beta} \right) \sum_{j=0}^b \frac{\nk {k \choose j} {n-k \choose k-j}}{2} \\
&=& \frac{1}{2}\nk p^4 \sum_{j=0}^b \kj \nkj + 2\nk p^3 { \left( \sum_{j=0}^b \kj \nkj \right) }^2.
\end{eqnarray*}
Clearly this quantity tends to zero when 
$$ p = o {\left( \frac{1}{\nk {\left( \sum_{j=0}^b \kj \nkj \right) }^2 } \right) }^{\frac{1}{3}}, $$
which was the stated restriction on $p$.
\hfill\qed

\medskip

\noindent {\bf Discussion.}  It has been shown in Theorem 5 that a multivariate Poisson approximation is valid for ${\cal L}(X_0,\ldots, X_b)$ when $p\ll\lr\nk \left( \sum_{j=0}^b \kj \nkj \right)^2 \rr^{-1/3}$.  This is a stronger condition, naturally, than those obtained in Theorem 4 for the univariate Poisson approximations for ${\cal L}(X_r)$.  We need to verify, however, that the threshold for multivariate Poisson approximation occurs at a level that is larger than the threshold for the EKR property; in other words we must have 
\begin{equation}\lr\nk{{n-k}\choose{k}}\rr^{-1/2}\le\lr\nk \left( \sum_{j=0}^b \kj \nkj \right)^2 \rr^{-1/3}.\end{equation}
The reason for imposing this requirement (we are not {\it required} to do so) is that one wishes to compute multivariate probability approximations for quantities involving {\it all} the $X_j$s; $0\le j\le b$, and, moreover, one would like to be able to meaningfully incorporate threshold situations into the multivariate approximation. In practice, given the value of $k$, (4) defines a condition that tells us how large $b$ may be before a Poisson approximation becomes unrealistic.  In general, the larger $k$ is, the larger $b$ is allowed to be.  This may be best seen by reexpressing (4) as
$$\nk^{1/2}{{n-k}\choose{k}}^{3/2}\ge\lr\sum\kj\nkj\rr^2.$$
\medskip

\noindent{\bf Acknowledgment}  The research of the first and third named authors was  supported by NSF Grants DMS-9619889.  The research of the second named author was supported by NSF Grants DMS-9619889 and DMS-0139291.


\medskip

\begin{thebibliography}{99}
\bibitem{alon} N.~Alon and J.~Spencer, {\it The Probabilistic Method, 2nd Edition,} John Wiley, New York, 2000.
\bibitem{bhj} A.~Barbour, L.~Holst and S.~Janson, {\it Poisson Approximation,} Oxford University Press, 1992.
\bibitem{ekr} P.~Erd\H os, C.~Ko, and R.~Rado (1961).  ``Extremal problems among subsets of a set," {\it Quart. J. Math. Oxford, Ser. 2} {\bf 12}, 313--318.
\bibitem{janson} A.~Godbole and S.~Janson (1996).  ``Random covering designs," {\it Journal of Combinatorial Theory, Series A} {\bf 75}, 85--98. 
\bibitem{greene} J.~Greene (2002). ``A new short proof of Kneser's conjecture," 
{\it Amer. Math. Monthly} {\bf 109}, 918--920.
\bibitem{wilson} J.~van Lint and R.~Wilson, {\it A Course in Combinatorics,} Cambridge University Press, 1992.
\end{thebibliography}
\end{document}